\documentclass{article}
\usepackage[utf8]{inputenc}
\usepackage[T1]{fontenc}
\newcommand\ignore[1]{}
\usepackage{fullpage} 
\usepackage{parskip}  
\usepackage{tabularx} 

\usepackage{amsmath,amssymb, amsfonts, amstext}
\usepackage{optidef}

\usepackage{commath} 

\usepackage{nth}
\usepackage{nicefrac}
\usepackage{siunitx}
\DeclareSIUnit\feet{ft}
\DeclareSIUnit\calorie{cal}

\usepackage{authblk}

\usepackage{booktabs}
\usepackage{tabularx}
\usepackage{graphicx}
\usepackage{float}
\usepackage{subfig}
\usepackage{wrapfig}
\usepackage{nth}
\usepackage{appendix}

\usepackage{xcolor}
\usepackage{subfiles}

\usepackage{todonotes}

 \usepackage[symbol]{footmisc}
 
\newcommand{\ie}{\emph{i.e., }}
\newcommand{\eg}{\emph{e.g., }}

\DeclareMathOperator{\Ebb}{\mathbb{E}}
\DeclareMathOperator{\Rbb}{\mathbb{R}}

\newcommand{\E}[1]{\Ebb\sbr{ #1 }}   
\newcommand{\Egiven}[2]{\Ebb\sbr {\giventhat{#1}{#2} }}
\newcommand{\Espace}[3]{\Ebb_{#1} \sbr {\giventhat{#2}{#3}}}

\newcommand{\inner}[2]{\left< #1, #2 \right>}
\newcommand{\giventhat}[2]{#1\;\middle|\;#2}
\newcommand{\bs}[1]{\boldsymbol{#1}}

\newcommand*\from{\colon}
\newcommand{\map}[3]{#1\from #2 \to #3}
\newcommand{\mapself}[2]{\map{#1}{#2}{#2}}

\newcommand{\code}[1]{$\mathtt{#1}$}

\title{The Koopman Expectation: An Operator Theoretic Method for Efficient Analysis and Optimization of Uncertain Hybrid Dynamical Systems}

\author{Adam R.~Gerlach}
\affil{United State Air Force Research Laboratory, Wright–Patterson Air Force Base, OH}
\author{Andrew Leonard}
\author{Jonathan Rogers}
\affil{Georgia Institute of Technology, Atlanta, GA}
\author{Chris Rackauckas}
\affil{Massachusetts Institute of Technology, Cambridge, MA}

\date{}
\begin{document}

\maketitle
\begin{abstract}
For dynamical systems involving decision making, the success of the system greatly depends on its ability to make good decisions with incomplete and uncertain information. By leveraging the Koopman operator and its adjoint property, we introduce the Koopman Expectation, an efficient method for computing expectations as propagated through a dynamical system. Unlike other Koopman operator-based approaches in the literature, this is possible without an explicit representation of the Koopman operator. Furthermore, the efficiencies enabled by the Koopman Expectation are leveraged for optimization under uncertainty when expected losses and constraints are considered. We show how the Koopman Expectation is applicable to discrete, continuous, and hybrid non-linear systems driven by process noise with non-Gaussian initial condition and parametric uncertainties. We finish by demonstrating a 1700x acceleration for calculating probabilistic quantities of a hybrid dynamical system over the naive Monte Carlo approach with many orders of magnitudes improvement in accuracy.

\end{abstract}

\section{Introduction \& Motivation}
The modern study of dynamics is dominated by the \emph{state space} view as pioneered by Poincaré \cite{goroff_new_1992} in the early 19th century \cite{nolte_tangled_2010}. By focusing on the \emph{dynamics of states}, this view offers intuitive geometric tools for studying dynamic systems and its impact on the sciences cannot be understated. In fact, some consider the state space view as ``one of the most powerful inventions of modern science'' \cite{gleick_chaos_2008}. However, despite its natural applicability to large classes of problems, it is ill-suited for high-dimensional \cite{andronov_theory_2011}, ill-described\footnote{For example, systems without explicit differential equations, \ie data-driven models}\cite{jones_whither_2001}, and uncertain systems \cite{budisic_applied_2012}. 

An alternative to the state space view of Poincaré is the \emph{observable space} view introduced by Koopman \cite{koopman_hamiltonian_1931} in 1931. Poincaré's view considers the evolution of points in the state space, while Koopman's view considers the evolution of observables, or functions, of the state space \cite{budisic_applied_2012}. This observable space view directly inspired Neumann's mean ergodic theorem \cite{halmos_legend_1973} and further led to the notion of the \emph{spectrum of a dynamical system} \cite{koopman_dynamical_1932} in 1932. In particular, the spectral properties of the so-called Koopman operator of a system reveals the mechanical properties of the underlying system itself \cite{koopman_dynamical_1932}. The Koopman operator is a linear operator on the space of scalar-valued functions of the system states. It describes the evolution of these functions as they are driven by the dynamics of the underlying system \cite{narasingam_koopman_2019}. So, the observable space view enables linear analysis of non-linear systems, albeit at the expense of dealing with potentially infinite-dimensional (function) spaces. 

Despite its initial scientific impact, this observable space view and the resulting Koopman operator had limited practical impact until its revival by the work of Mezic \cite{mezic_spectral_2005-2, mezic_comparison_2004} in the early 2000s. This ignited an outburst in theoretical advancements and computational tools for analyzing, designing, and controlling dynamic systems via the Koopman operator. It is easy to understand why, as the Koopman operator promises to:

\begin{enumerate}
	\item capture full nonlinear dynamics in a linear setting \cite{budisic_applied_2012}
	\item not sacrifice information, like with traditional linearization techniques \cite{budisic_applied_2012}
	\item provide a bridge between nonlinear systems and existing linear tools for analysis, design, and control \cite{brunton_koopman_2016,kutz_dynamic_2016}
\end{enumerate} 

Despite these promises, to date, the Koopman operator has primarily found practical application in the context of data-driven methods for both ill- and well-defined\footnote{In the case of well-defined systems, analytical models are used to generate synthetic datasets and data-driven methods are then employed on these datasets} systems. This is due to the fact that in order to study the spectral properties of the Koopman operator for a given system, a representation of the operator is required. Except for the simplest of systems, the Koopman operator is typically not known explicitly.  However, simple numerical methods exist for computing reduced order approximations of the operator directly from data \cite{chen_variants_2012, korda_data-driven_2020, kutz_dynamic_2016, schmid_dynamic_2010, tu_dynamic_2014, williams_datadriven_2015}.

In contrast to this approach, we are interested, not in the Koopman operator directly, but the \emph{action of the operator} on observables of interest. We intend to focus this in the context of optimization under uncertainty by exploiting one particular property of the Koopman operator: its adjoint property. The Koopman operator (the \emph{pull-back} operator on observables) is the adjoint to the Frobenius-Perron (FP) operator (the \emph{push-forward} operator on densities), \ie
\begin{equation*}
	\inner{Pf}{g} = \inner{f}{Ug}
\end{equation*}
where $\inner{\cdot}{\cdot}$ is the inner product, $P$ is the FP (push-forward) operator applied to the density $f$, and $U$ is the Koopman (pull-back) operator applied to the observable of interest $g$. 

Although this property is often cited in the literature, it was never\footnote{To the best of the authors' knowledge} used outside of making mathematical proofs about the Frobenius-Perron operator \cite{lasota_chaos_2013} until the recent work of the authors \cite{leonard_probabilistic_2019-1,leonard_koopman_2019-1,leonard_probabilistic_2020-1}. When considering probability density functions, one can reinterpret this adjoint relationship in terms of expectations:
\begin{equation*}
	\Egiven{g(\bs X)}{\bs X\sim Pf} = \Egiven{Ug(\bs X)}{\bs X \sim f}
\end{equation*}
In other words, to compute the expected value of some quantity of interest, $g$, of an uncertain system in the future, one could first push the uncertainty forward through the system $\del{Pf}$ and then compute the expectation of $g$ (the LHS, \emph{Frobenius-Perron Expectation}). Alternatively, one could pull the same quantity of interest supported at some future time back to the initial time via the Koopman operator leading to a new quantity of interest $\del{Ug}$. The expectation of this new quantity of interest is then computed with respect to the initial density $\del{f}$ (the RHS, \emph{Koopman Expectation}). Although these provide different means to the same end, the Koopman Expectation (the RHS) provides significant numerical advantages \cite{leonard_probabilistic_2019-1,meyers_koopman_2019-1}. When used in the context of optimization under uncertainty, the optimization problem can be written in the so-called \emph{Koopman Expectation Form} to exploit these advantages. 

In contrast to other Koopman operator-based approaches in the literature, solving the Koopman Expectation  is possible without explicit knowledge of the  Koopman operator $U$ itself. This stems directly from the definition of the Koopman operator: 
\begin{equation*}
	Ug\del{\bs x} = g\del{S\del{\bs x}}
\end{equation*}
where $S$ maps the system states to some instance in the future \cite{lasota_chaos_2013}. Here, an explicit representation of the Koopman operator is not necessary. Hence, we say that the \emph{action of the operator} is known for a given function $g$. To further understand the distinction between the action of the operator vs.~the operator itself, compare this with the familiar derivative operator. One can easily apply the derivative operator to a given function without knowing the operator's spectral properties, \ie the action of the operator is known for the given function.

In the following section, we provide the reader with the technical background for understanding the Koopman Expectation and its benefits in a general context. 

\section{Background: The Frobenius-Perron \& Koopman Operators}
To motivate the application of the Koopman operator for optimization of uncertain systems, we first introduce the Frobenius-Perron (FP) operator, a linear transfer operator on densities \cite{lasota_chaos_2013}. Let $\mapself{S}{\Rbb^n}$ be a dynamical system. Common choices are where the underlying dynamical process is a discrete process
\begin{equation}
  y(n) = f(y(n-1))
\end{equation}
or a continuous dynamical system described by an ordinary differential equation
\begin{equation}
    \dot{y} = f(y)
\end{equation}
For these equations, we define $S(x) = y(T)$ where $y(0) = x$. We wish to compute probabilistic statements on functions, or observables, of the dynamical system with respect to uncertainty of the initial condition. We note in passing that uncertainty with respect to parameters can be incorporated into this framework by extending the dynamical system with states having null dynamics and initial conditions equal to the associated parameters. 

\begin{figure}[h]
	\centering
	\includegraphics[width=0.8\textwidth]{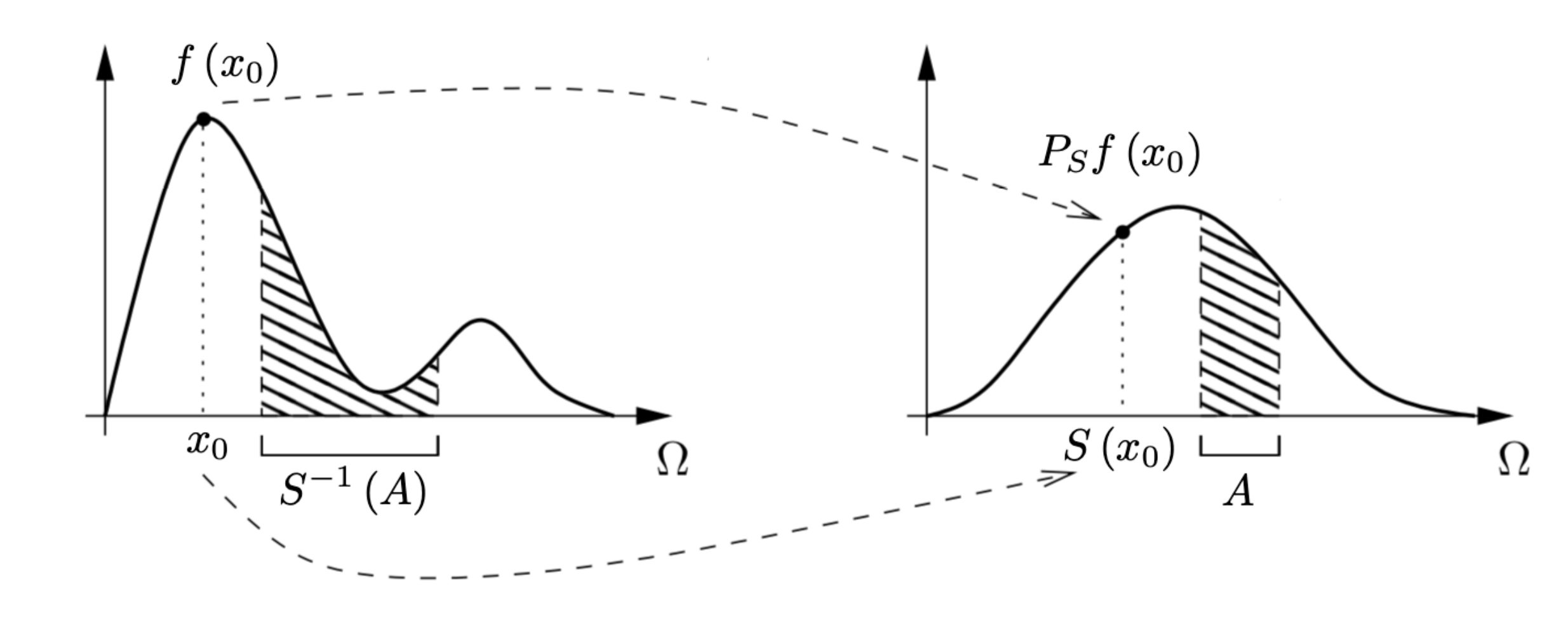}
	\caption{Graphical representation of Eq.~\ref{eq:FP}, where shaded regions have equal area. (Figure adapted from \cite{weise_global_2009})}
	\label{fig:FP_transport}
\end{figure}

Without loss of generality, we model the uncertainty of the initial condition via the standard probability space $\del{\Omega, \mathcal{A}, \mu}$ where $\Omega = \Rbb^n$. For a given dynamical system $\mapself{S}{\Rbb^n}$, its associated FP  operator, $P_S$, is defined such that the following equality is satisfied: 
\begin{equation}
\int_A P_S f \del{\bs x}\mu\del{\dif\bs x} = \int_{S^{-1}\del{A}}f\del{\bs x} \mu\del{\dif{\bs x}},\quad \forall A\in\mathcal{A} \label{eq:FP}
\end{equation}
where $S^{-1}\del{A}$ is the counter-image of $A$ and $f$ is a probability density function \cite{lasota_chaos_2013}. If $S$ is both measurable and nonsingular, then $P_S$ is uniquely defined by Eq.~\ref{eq:FP} \cite{lasota_chaos_2013}. Intuitively, this means that there exists a unique linear operator $P_S$ such that the probability mass of subregion $S^{-1}\del{A}$ is conserved when this subregion is transported to subregion $A$ via $S$. This is true for all possible subregions of $\Omega$. This is illustrated in \figref{fig:FP_transport}. 

When $S$ is both differentiable and invertible, then the transformed density in Eq.~\ref{eq:FP} can be solved for explicitly as
\begin{equation}
P_S f\del{\bs x} = f\del{S^{-1}\del{\bs x}}\abs{\od{S^{-1}\del{\bs x}}{\bs x}} \label{eq:FP_solved}
\end{equation}
where $\abs{\nicefrac{\dif S^{-1}\del{\bs x}}{\dif \bs x}}$ is the determinant of the Jacobian of $S^{-1}$ \cite{lasota_chaos_2013}. This indicates that the transformed density of $\bs x$ is the initial density at its counter-image scaled by the factor $\abs{\nicefrac{\dif S^{-1}\del{\bs x}}{\dif \bs x}}$. This factor enforces the conservation of probability mass by accounting for the local expansions and contractions of the space due to $S$ \cite{weise_global_2009}. 





 Let us now turn our attention to the Koopman operator. Whereas the FP operator provides a mechanism to transport densities through maps or dynamical systems, the Koopman operator provides a mechanism to transport functions \cite{koopman_hamiltonian_1931}. For a non-singular map $S$, the Koopman operator, $U_S$ is defined as the composition of an observable $g$ and $S$, \ie
\begin{equation}
U_S g\del{\bs x}=g\del{S\del{\bs x}},\quad \forall g\in L^\infty \label{eq:Koopman}
\end{equation}

%
%

\section{The Koopman Expectation} \label{sec:koopman_expectation}
For a map $S$, the corresponding Koopman operator, $U_S$, possesses a particular property of interest: it is adjoint to the FP operator, $P_S$, \ie 
\begin{equation}
\inner{P_Sf}{g}=\inner{f}{U_S g},\quad \forall f\in L^1, \forall g\in L^\infty \label{eq:adjoint}
\end{equation}
where $\inner{\cdot}{\cdot}$ is the inner product. 

For real-valued functions, $g$, on $\Omega$, the inner products in Eq.~\ref{eq:adjoint} can be reinterpreted as  expectations, \ie
\begin{equation} \label{eq:inner_prod_integral}
\Egiven{v\del{\bs X}}{\bs X\sim u}=\int_\Omega u\del{\bs x}v\del{\bs x}\dif\bs x=\inner{u}{v}
\end{equation}
where the measure $\mu$ is induced by the density $u$ \footnote{\ie $\int \mu\del{\dif \bs x}=\int u\del{\bs x}\dif \bs x$}.Thus, Eq.~\ref{eq:adjoint} can be rewritten as
\begin{equation}
\Egiven{g\del{\bs X}}{\bs X\sim P_Sf} = \Egiven{U_S g\del{\bs X}}{\bs X \sim f}
\end{equation}
or for notational simplicity as
\begin{equation}
\Egiven{g}{ P_Sf} = \Egiven{U_S g}{f} \label{eq:Exp_def}
\end{equation}
We refer to the left- and right-hand sides of Eq.~\ref{eq:Exp_def} as the \emph{Frobenius-Perron} (FP) and \emph{Koopman Expectations}, respectively. 

\begin{figure}[h]
	\centering
	\includegraphics[width=0.5\linewidth]{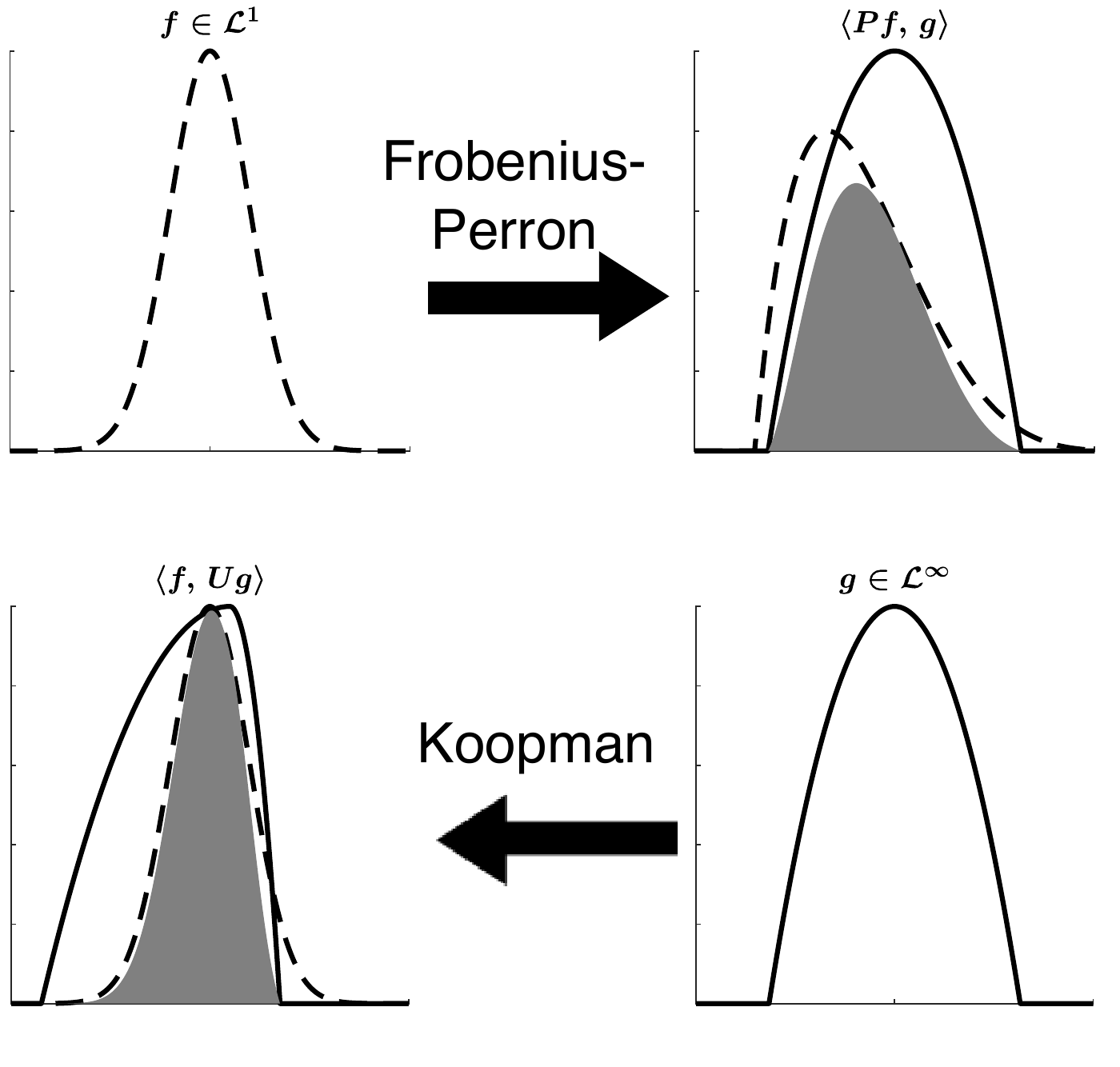}
	\caption{Illustration of the FP and Koopman operator adjoint property. The inner products, represented by the area of the filled regions, are equivalent. (Figure reproduced from \cite{leonard_probabilistic_2019-1})}
	\label{fig:kfpdrawing}
\end{figure}

\figref{fig:kfpdrawing} provides a 1D illustration of Eq.~\ref{eq:Exp_def}. The top row represents the FP Expectation while the bottom row represents the Koopman Expectation. For the top row, the PDF $f$ (dashed line) is pushed to the right through the system dynamics via $P_S$ and an inner product is taken with $g$ (solid line). The expected value $\Egiven{g}{P_Sf}$ is represented by the area of the shaded region. Conversely, on the bottom row, the function $g$ is pulled to the left through the system dynamics via $U_S$. The expected value $\Egiven{U_Sg}{f}$ is represented by the area of the shaded region. The areas of the two shaded regions are equal. In summary, each form provides a different means to the same end, the expectation. 

For problems where an expectation is required, the Koopman Expectation (right-hand side of Eq.~\ref{eq:Exp_def}) offers numerous computational advantages:

\begin{enumerate}
	\item The simple evaluation of  $U_Sg$ vs.~$P_Sf$ (Eqs.~\ref{eq:Koopman} and \ref{eq:FP_solved}, respectively). $U_Sg$ is the direct composition of the function $g$ and map $S$. In contrast, computing $P_Sf$ requires knowing the inverse map $S^{-1}$, which is oftentimes unknown. This is especially true when $S$ represents a cyber-physical system or is evaluated as a black-box. 
	\item For asymptotically stable systems, the support of the transformed PDF collapses to a manifold. Thus, the volume of this support approaches $0$ over time. Consequently, in order to conserve probability mass, $\sup P_Sf\rightarrow\infty$, leading directly to numerical overflow in computer implementations \cite{hoogendoorn_uncertainty_2018}. 
	\item For systems with uncertainty in only a subset of the initial conditions, the push-forward operation can lead to drift across the entire state-space. This requires the FP Expectation to be taken over the entirety of the state-space, whereas the Koopman Expectation can be taken over only the subset of the state-space with initial uncertainty. 
	\item Because $U_Sg$ is supported on the original domain, the structure of the domain is preserved. One simple consequence of this is that the bounds of integration are well defined for the Koopman Expectation, unlike for the FP Expectation. This is a documented issue with FP in the literature \cite{probe_new_2016}. Secondly, by preserving the structure of the domain, structures that lead to simpler solution approaches can be prescribed, \eg quadrature integration.
	\item If the expectation of multiple functions with varying supports in space-time is required, the pull-back of each to a common domain results in a single, vector-valued expectation calculation. 
\end{enumerate}

\section{Higher-Order Statistics via the Koopman Expectation} \label{sec:statistics}
The power of the Koopman operator is enabled by the generality of the observables, or space thereof, that it operates on. However, when coupled with uncertainty, the methodology above only allows for analysis in the \textit{mean} sense instead of through richer information about the transformed uncertainty; as provided by the FP operator and other push-forward methods. Nevertheless, through the specific selection of a set of observables, similar information may be reconstructed via the Koopman Expectation. The process for determining the required observables is demonstrated by example for central moments, covariances, and correlations below. This process can be readily extended to additional higher-order statistics that can be decomposed into a basis of expectations.

\subsection{Central Moments} \label{sec:central_moments}
The $n^{th}$ central moment, $\hat{m}_n$, is defined as
\begin{equation} \label{eq:central_moment_n}
\hat{m}_n = \E{ \left(Z - \mu(Z) \right)^n }
\end{equation}
where $\mu(Z) = \E{Z}$ is the mean, or first raw moment, of the random variable $Z$. 

By expanding the interior binomial, Eq.~\ref{eq:central_moment_n} can be written as
\begin{equation}
\hat{m}_n = \E{ \sum_{k = 0}^n \binom{n}{k}(-1)^{n-k} \: \mu(Z)^{n-k} Z^k }
\end{equation}
Due to the linearity of the expected value operator and the fact that the mean is a constant, this can be rewritten as
\begin{equation} \label{eq:central_moment_expanded}
\hat{m}_n = \sum_{k = 0}^n \binom{n}{k}(-1)^{n-k} \mu(Z)^{n-k} \: \E{Z^k} 
\end{equation}

Equation \ref{eq:central_moment_expanded} shows that to compute the $n^{\mathrm{th}}$ central moment, the expected values of the $n$-powers of the random variable are needed. This defines the required set of $n$ observable functions to be use with the Koopman Expectation as
\begin{equation}
\left\{z, \: z^2, \: z^3, \, \dots, \: z^n \right\}
\end{equation}
\subsection{Covariance and Correlation} \label{sec:covariance_correlation}
In a similar fashion, the covariance of two observables may be computed using the Koopman Expectation. The covariance of the random variables $Z_1$ and $Z_2$ is defined as
\begin{equation} \label{eq:covariance}
\mathrm{cov} \left( Z_1, \, Z_2 \right) =  \E{Z_1 Z_2} - \E{Z_1} \E{Z_2}
\end{equation}
It follows directly that the set of observable functions required are
\begin{equation*}
\left\{z_1, \:z_2, \: z_1 z_2 \right\}
\end{equation*}

Additionally, the correlation between two random variables is defined as
\begin{equation} \label{eq:correlation}
\rho_{Z_1, \, Z_2} = \frac{\mathrm{cov}(Z_1, \, Z_2)}{\sigma_{Z_1} \sigma_{Z_2}} = \frac{\E{Z_1 Z_2} - \E{Z_1} \E{Z_2}}{\sqrt{\E{Z_1^2} - \E{Z_1}^2 } \sqrt{\E{Z_2^2} - \E{Z_2}^2 }}
\end{equation}
where $\sigma_Z$, the standard deviation of random variable $Z$, is defined as the square root of the variance, or second central moment. It follows directly form the right hand side of Eq.~\ref{eq:correlation} that the set of observable functions required are
\begin{equation*}
\left\{z_1, \: z_1^2, \: z_2, \:z_2^2, \: z_1 z_2 \right\}
\end{equation*}

\section{Incorporating Process Noise} \label{sec:process_noise}

Our previous applications of the Koopman Expectation did not include process noise. We can expand the framework to include process noise by utilizing Koopman operator theory on random dynamical systems \cite{vcrnjaric2019koopman}. A random dynamical system on a measurable space $(X,\mathcal{B})$ over $(\mathcal{\theta}_t)_{t\in T}$ on $(\Omega,\mathcal{A},\mu)$ is a measurable map
\begin{equation}
    S: T \times X \times \Omega \rightarrow X,
\end{equation}
such that $S(0,x_0,\omega) = x_0$ and
\begin{equation}
    S(t+s,x_0,\omega) = S(s,S(t,x_0,\omega),\mathcal{\theta}_s \omega),
\end{equation}
for all $t,s \in T$ and for all $\omega$ outside a $\mu$-nullset \cite{arnold1995random} (for simplicity of notation, we refer to the starting time as $t=0$). A common case of a random dynamical system is a random ordinary differential equation (RODE)
\begin{equation}
    \dot{y}(t) = F(\mathcal{\theta}_t \omega,y(t)).
\end{equation}
where in less technical terms, $\mathcal{\theta}_t \omega$ is a random draw of a state-independent forcing process. In this case, the same formulation of the Koopman Expectation carries over where $S(t,x,\omega) = S(t,x)$ corresponds to deterministic dynamics. Thus the extension of the Koopman Expectation on a random dynamical system is simply to include the parameterization of the forcing process as part of the integrated probability space.

We will illustrate the process of extending the Koopman Expectation to the common case of Gaussian white noise, \ie where $\mathcal{\theta}_t \omega$ represents a Wiener process under the Ito definition, but the discussion similarly applies to other noise processes. An Ito-driven RODE can be rewritten as a stochastic differential equation (SDE):
\begin{equation}
    \dif y(t) = \phi(y(t))\dif t + \psi(y(t))\dif W_t
\end{equation}
where $W_t$ is the random forcing process known as the Wiener process \cite{han2017random}. The forcing process $W_t$ can be approximated in multiple ways. If the dynamical system is approximated by a fixed time stepping process, then we can represent $W_t$ exactly by the summation:
\begin{equation}
  W(n \Delta t) = \sum_{i=1}^{n} \Delta W_i
\end{equation}
where $\Delta W_i \sim N(0,\Delta t)$. If a continuous $W(t)$ is needed, for example when approximating the solution to the SDE by adaptive time stepping \cite{rackauckas2017adaptive}, then this process can be linearly interpolated between the defined points as justified by the fact that the mean of the Brownian bridge process is the linear form. An alternative formulation of a continuous Brownian motion via finite random variables is to use the fact that $W(t)$ on $[0,T]$ can be represented by a truncation of the Karhunen–Loève expansion \cite{stark1986probability,dutta_nonlinear_2015}:
\begin{equation}
    W_t = \sqrt{2} \sum_{k=1}^\infty Z_k \frac{\sin \left(\left(k - \frac 1 2 \right) \pi \frac{t}{T}\right)}{ \left(k - \frac 1 2 \right) \pi}
\end{equation}
where $Z_k$ are independent and identically distributed standard normal variables. Using this formulation we can represent the Koopman Expectation over a finite dimensional probability space and thus apply the previous methods.

An alternative formulation for the incorporating of the process noise is to utilize the expectation operator to define a deterministic dynamical system from the random dynamical system \cite{lasota_chaos_2013}. If one defines:
\begin{equation}
    \tilde{S} = \mathbb{E}[S]
\end{equation}
then one can show that the resulting map $\tilde{S}$ is a deterministic dynamical system of which the previous theory can be directly applied to. We note in passing that this formulation may be more computationally expensive as it amounts to computing the solution to many SDE systems per integrand calculation, whereas the direct discretization only requires one.

\section{Optimization Under Uncertainty via the Koopman Expectation}
\label{sec:koop_opt}

Equipped with knowledge of the Koopman Expectation, we now discuss how it can be used for optimization under uncertainty. 

Define the decision variables as $\bs u \in \mathcal{U}$, where $\mathcal{U}$ is the set of admissible decisions. Next, consider the initial joint PDF of the system conditioned on $\bs u$, $\map{f_0\del{\giventhat{\cdot}{\bs u}}}{\Omega}{\Rbb^+}$, where the form of $f_0\del{\giventhat{\cdot}{\bs u}}$ is problem specific. We define an \emph{Initial Condition Space}, $\mathcal{I}$, as the support of $f_0\del{\giventhat{\cdot}{\bs u}}$, \ie
\begin{equation}
\mathcal{I}=\set{\bs x \from f_0\del{\giventhat{\bs x}{\bs u}}>0, \bs x \in \Omega, \bs u \in \mathcal{U}} \label{eq:ic_space}
\end{equation}
 Next, consider an objective function $\map{g}{\mathcal{O}}{\Rbb}$. We call $\mathcal{O}$ the \emph{Objective Space} and define it implicitly as
\begin{equation}
\mathcal{O}=\set{ \bs x \from h\del{\bs x} = 0, \bs x \in \Omega} \label{eq:obj_space}
\end{equation}
where $\map{h}{\Omega}{\Rbb}$ is problem specific. Similarly, consider a series of $m$ constraint functions $\map{c_i}{\mathcal{C}_i}{\Rbb}$. We call $\mathcal{C}_i$ the \emph{$i$-th Constraint Space} and define it implicitly as
\begin{equation}
\mathcal{C}_i=\set{\bs x \from q_i\del{\bs x}=0, \bs x \in \Omega}, \quad i=1,\ldots,m
\end{equation}
where $\map{q_i}{\Omega}{\Rbb}$ is also problem specific. It is important to note that in most applications of interest these spaces differ. 

We then formulate the optimization problem as 
\begin{argmini}|s|
	{\bs u \in \mathcal{U}}{\Ebb_\mathcal{O}\sbr{g\del{\bs X}}}{}{\bs u^*=}
	\addConstraint{\Ebb_{\mathcal{C}_i}\sbr{c_i\del{\bs X}}}{\leq \lambda_i,\quad}{i=1,\ldots,m} 
\end{argmini}

We introduce the subscript notation $\Ebb_\mathcal{O}$ to make explicit the space over which the expectation is taken. When $\mathcal{O}, \mathcal{C}_i\nsubseteq\mathcal{I}$, the push-forward of $f_0$ is required to evaluate these expectations. If we consider the mappings defined by $\map{S}{\Omega}{\Omega}$ such that $S\del{\mathcal{I}}\subseteq\mathcal{O}$ and $\map{Q_i}{\Omega}{\Omega}$ such that $Q_i\del{\mathcal{I}}\subseteq\mathcal{C}_i$ we can be more explicit with our formulation and rewrite it in the FP Expectation Form, \ie
\begin{argmini}|s|
	{\bs u \in \mathcal{U}}{\Espace{\mathcal{O}}{g\del{\bs X}}{\bs X\sim P_Sf_0\del{\giventhat{\cdot}{\bs u}}}}{}{\bs u^*=}
	\addConstraint{\Espace{\mathcal{C}_i}{c_i\del{\bs X}}{\bs X\sim P_{Q_i}f_0\del{\giventhat{\cdot}{\bs u}}}}{\leq \lambda_i,\quad}{i=1,\ldots,m}  \label{eq:opt_fp_form}
\end{argmini}
where $P_S$ and $P_{Q_i}$ are the FP operators for $S$ and $Q_i$, respectively. Note that $S$ and $Q_i$ are defined such that they relate the Initial Condition Space to the Objective and Constraint Spaces. Furthermore, this relates the set of admissible decisions to the Objective and Constraint spaces via their conditional relationship with the Initial Condition Space (Eq.~\ref{eq:ic_space}).

Using the adjoint relationship, this is rewritten in the Koopman Expectation Form as
\begin{argmini}|s|
	{\bs u \in \mathcal{U}}{\Espace{\mathcal{I}}{U_Sg\del{\bs X}}{\bs X\sim f_0\del{\giventhat{\cdot}{\bs u}}}}{}{\bs u^*=}
	\addConstraint{\Espace{\mathcal{I}}{U_{Q_i}c_i\del{\bs X}}{\bs X\sim f_0\del{\giventhat{\cdot}{\bs u}}}}{\leq \lambda_i,\quad}{i=1,\ldots,m}  \label{eq:opt_koop_form}
\end{argmini}

where $U_S$ and $U_{Q_i}$ are the Koopman operators for $S$ and $Q_i$, respectively. Note that in the Koopman Expectation Form all $m+1$ expectations are taken in the Initial Condition Space. This is unlike the FP Expectation Form, where all expectations are generally taken in separate spaces. Further note that $f_0\del{\giventhat{\cdot}{\bs u }}$ is used directly in the Koopman Expectation Form, whereas $m+1$ push-forwards of $f_0\del{\giventhat{\cdot}{\bs u }}$ are required in the FP Expectation Form.

\section{Computational Considerations} \label{sec:computational_details}
The optimization of dynamical systems is already a computationally intensive problem requiring high-performance differential equation solvers and adjoint methods for fast gradients. Optimization under uncertainty adds another layer to the computational complexity since we now need to evaluate and differentiate a multidimensional integral defined potentially over numerical differential equation solves. However, despite these lofty requirements, this problem can be made tractable by decomposing the Koopman Expectation and optimization into functional components that use just a handful of generally applicable computational tools. 

At the lowest level is the observable function \code{y = g(x)} which computes a scalar value \code{y} from the system state \code{x}. The system of interest is represented by its map, \code{x' = S(x)}, which computes the new system state \code{x'} from \code{x}. For discrete or continuous systems with closed form solutions, the map may be represented exactly. For more complicated systems, the function may, for example, hold iterative applications of a discrete system or a call to a numerical integrator for continuous systems. Independent of the choice of how \code{S} is implemented, its result is composed directly as the argument to \code{g}, \ie 
\begin{equation*}
	\mathtt{y = g\left( S(x) \right)}
\end{equation*}
For the chosen implementation of \code{S}, this composition represents the action of the Koopman operator on the function \code{g} \textit{exactly}. As such, there is no error introduced by the methodology itself. The only error introduced is through the user's chosen representations of the system map and observable function.

To compute the Koopman Expectation from Eq.~\ref{eq:inner_prod_integral}, the result from the above composition can be embedded into the integrand alongside a representation of the PDF of the initial condition, \code{f(x)}. In this form, the integrand may be point-wise evaluated as \code{I(x) = g(S(x)) * f(x)} and can be passed to a numerical integration method chosen by the user. 

\begin{figure}[h]
	\centering
	\includegraphics[width=0.85\linewidth]{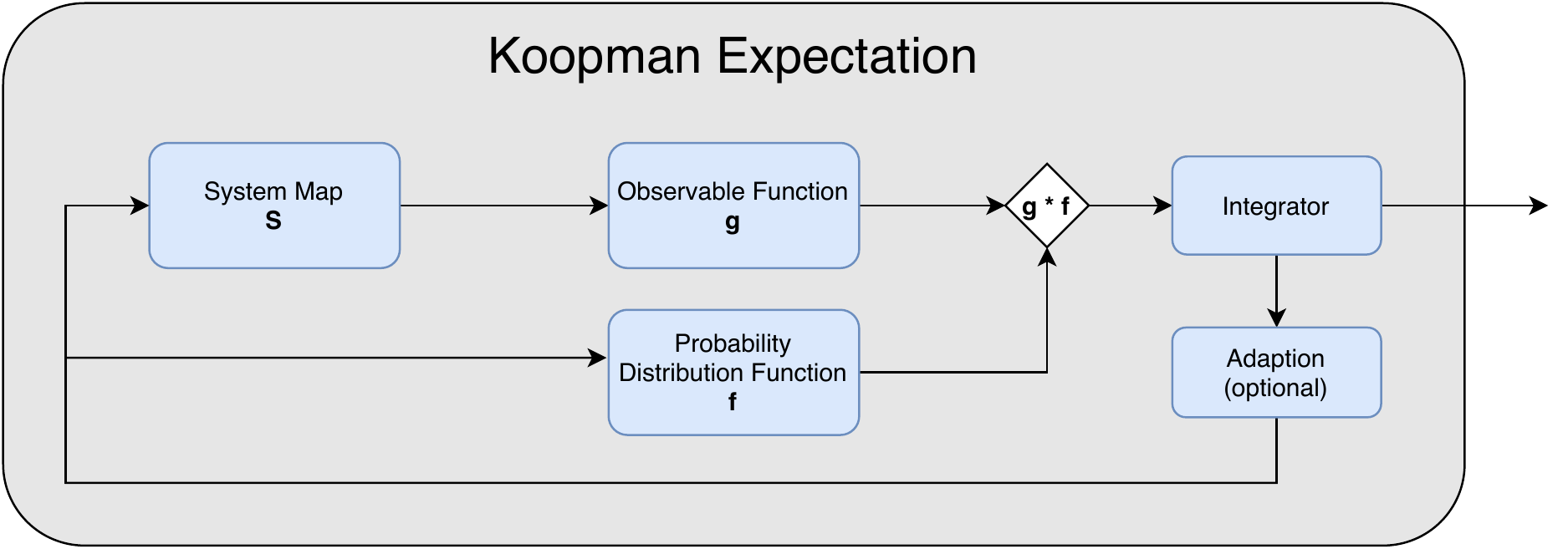}
	\caption{Functional decomposition of the Koopman Expectation}
	\label{fig:opt_flowchart}
\end{figure}

From these simple components, the full Koopman Expectation may be written functionally as 
\begin{equation*}
	\mathtt{sol = koopman\_expectation(S, g, f, integration\_method)}
\end{equation*}
where \code{S}, \code{g}, \code{f} are all functions of \code{x} and some possible model parameters and \code{sol} is the scalar-valued expectation. Figure~\ref{fig:opt_flowchart} illustrates this decomposition.

In the context of optimization under uncertainty, the \code{ koopman\_expectation} function can be further composed with a user-specified optimization method. For gradient-based methods the gradient of \code{ koopman\_expectation} can be computed by utilizing automatic differentiation over the quantity of interest, where specializations on the differential equation solver allow for replacing the propagation through the differential equation and multidimensional quadrature code with analytically defined adjoints for these operators. As such, the expectation of a quantity of interest and its gradient can be optimally implemented from the following primitives:

\begin{itemize}
	\item An automatic differentiation tool which allows for defining optimized adjoints for specific functions
	\item Optimized implementations of adjoint expressions for differential equations
	\item Optimized implementations of adjoint expressions for spatial integration methods, such as multidimensional quadrature
	\item Easily parallelizable differential equation solver with support for accelerator hardware such as GPUs
\end{itemize}

This capability is available in the DiffEqUncertainty.jl\footnote{https://github.com/SciML/DiffEqUncertainty.jl} package for the Julia programming language.

\section{Illustrative Example}

To demonstrate the performance of the Koopman Expectation in terms of accuracy and computational efficiency as compared to the standard Monte Carlo approach consider a 2D bouncing ball with an uncertain coefficient of restitution\footnote{The coefficient of restitution is the ratio of relative velocities just before and just after collision } \ie

\begin{equation}
	\ddot{\bs{x}}=\begin{bmatrix}
	\ddot{x}\\ 
	\ddot{z}
	\end{bmatrix} 
	=  
	\begin{bmatrix}
	0\\
	-g
	\end{bmatrix},\quad x_0=\SI{2}{\meter}, \dot{x}_0=\SI{2}{\meter\per\second}, z_0=\SI{50}{\meter},\dot{z}_0=\SI{0}{\meter\per\second}
\end{equation}

where $x$ and $z$ are the horizontal and vertical position of the ball, respectively, and $g$ is the acceleration due to gravity on the Earth's surface. 

For time instances where $z=0$, the vertical velocity of the ball is modified according to 
\begin{equation}
	\dot{z}^+=-\alpha\dot{z}^-
\end{equation}
where  $\dot{z}^+$ and $\dot{z}^-$ are the ball's vertical velocity just prior and after impact, respectively, and $\alpha$ is the normal coefficient of restitution. Uncertainty in $\alpha$ is modeled as a truncated normal distribution on the interval $\sbr{0.84,1}$ with mean $\mu_\alpha=0.9$ and standard deviation $\sigma_\alpha=0.02$. We denote the PDF of this distribution by $f_\alpha$.

\begin{figure}[h]
	\centering
	\includegraphics[width=0.65\linewidth]{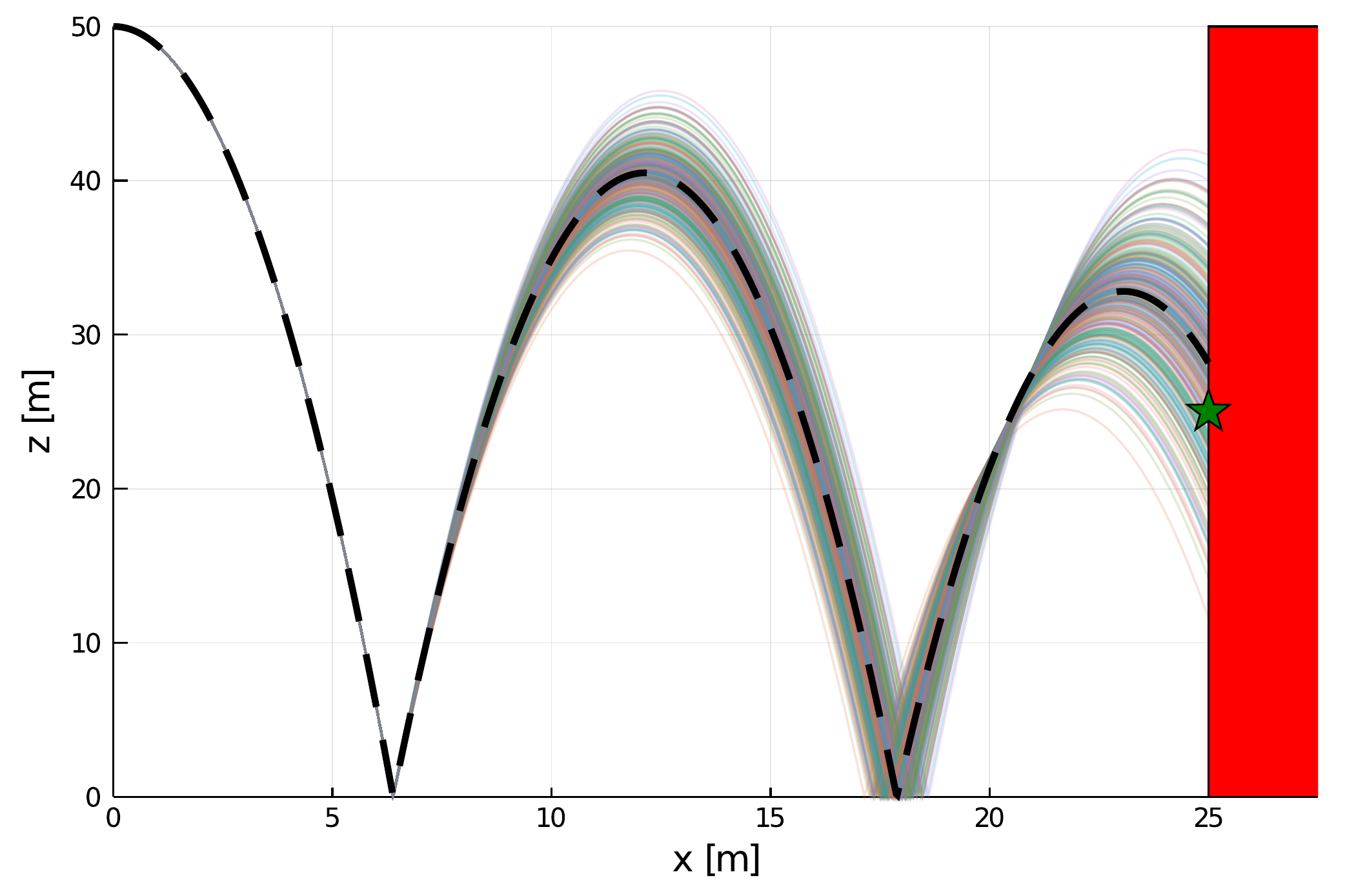}
	\caption{Bouncing Ball Scenario}
	\label{fig:ball_scenario}
\end{figure}

We wish to compute the expected squared miss distance from a target point on a vertical wall with coordinates $\del{x^*,z^*}=\del{\SI{25}{\meter}, \SI{25}{\meter}}$. This scenario is illustrated in Figure \ref{fig:ball_scenario} with the target point indicated by the green star. The nominal trajectory (dashed line) and 350 random trajectories are also shown. 

For this problem, $S$ maps the  initial condition to the vertical wall at $x=x^*$ as parameterized by $\alpha$.  $S$ is computed via numerical simulation using the Tsitouras 5/4 Runge-Kutta method \cite{tsitouras2011runge} provided by the DifferentialEquations.jl solver suite \cite{rackauckas_differentialequations.jl_2017} with event detection for both wall and ground impact. The simulation is terminated when the ball impacts the vertical wall. 

As we are interested in the squared miss distance from the target point on the wall, we define the observable $g$ as 
\begin{equation}
	g\del{\bs{x}}=\del{z-z^*}^2
\end{equation}
From the Koopman Expectation we have
\begin{align}
	\Egiven{g\del{S\del{\bs x,\alpha}}}{f_\alpha} = \int_{0.84}^{1}g\del{S\del{\bs x, \alpha}}f_\alpha\del{\alpha}d\alpha
\end{align}

The analytical solution to this expectation (Appendix \ref{ap:ball_analytical}) will serve as the truth value. For this particular example, we choose to use  h-adaptive integration \cite{genz_remarks_1980} to compute this integral in the Koopman Expectation. By doing so, the integration scheme adaptively selects values for $\alpha$ such that the relative and absolute error tolerences for the expectation are satisfied via an error estimate. Here, we set both of these tolerances to $\SI{1e-2}{}$. 

The efficiency of the Koopman Expectation is compared to the standard Monte Carlo simulation approach. All Monte Carlo simulations were performed by randomly sampling from $f_\alpha$ and evaluating $S$. Table \ref{tab:mc_v_koopman} shows the resulting expected value and computation time required from running \SI{100000}{} Monte Carlo simulations along with those from leveraging the Koopman Expectation. The residual computed by the h-adaptive integration is also reported for the Koopman Expectation. Each Monte Carlo simulation was run in parallel on a 6 core processor, whereas the Koopman Expectation was computed serially on the same processor. For this problem, the h-adaptive integration scheme used for the Koopman Expectation only required 15 simulations to reach the specified error tolerances, resulting in approximately a 1700x speed-up. Although ~6666x fewer simulations are required, a 6666x speed-up is not realized due to two primary factors:
\begin{enumerate}
	\item Overhead associated with the h-adaptive integration. This overhead is independent of the computational complexity of the system map $S$. So, for problems involving computationally expensive maps, the relative impact of this overhead diminishes. 
	\item Parallel versus serial implementations. As noted above, the Monte Carlo simulations were conducted in parallel on a 6 core processor, while the 15 simulations required for the Koopman Expectation were conducted serially. There is no theoretical or technical limitation preventing the batch parallel execution of these simulations. However, it does require a spatial integration library that supports batch integrand evaluations. Because of the overhead associated with parallel execution and integration and the limited number of simulations required to meet the integration tolerances of this problem, a parallel implementation offers minimal benefit here. However, for more complex problems parallel batch simulations can provide significant additional speed-ups, especially when coupled with massively parallel compute architectures like GPUs.
\end{enumerate}


\begin{table}[ht]
	\centering
	\caption{Bouncing Ball Monte Carlo Results}
	\label{tab:mc_v_koopman}
	\begin{tabular}[t]{lccc}
		\toprule
		& Analytical & Monte Carlo & Koopman \\ 
		\midrule
		Number of Simulations &- &\SI{100000}{} & 15 \\ 
		Expected Value $\sbr{\si{\meter\squared}}$&36.008 &35.782 & $36.008\pm \SI{4.987e-4}{}$\\ 
		Computation Time $\sbr{\si{\second}}$& -&2.060 & 0.0012 \\ 
		\bottomrule
	\end{tabular}
\end{table} 

\begin{figure}[h]
	\centering
	\includegraphics[width=0.65\linewidth]{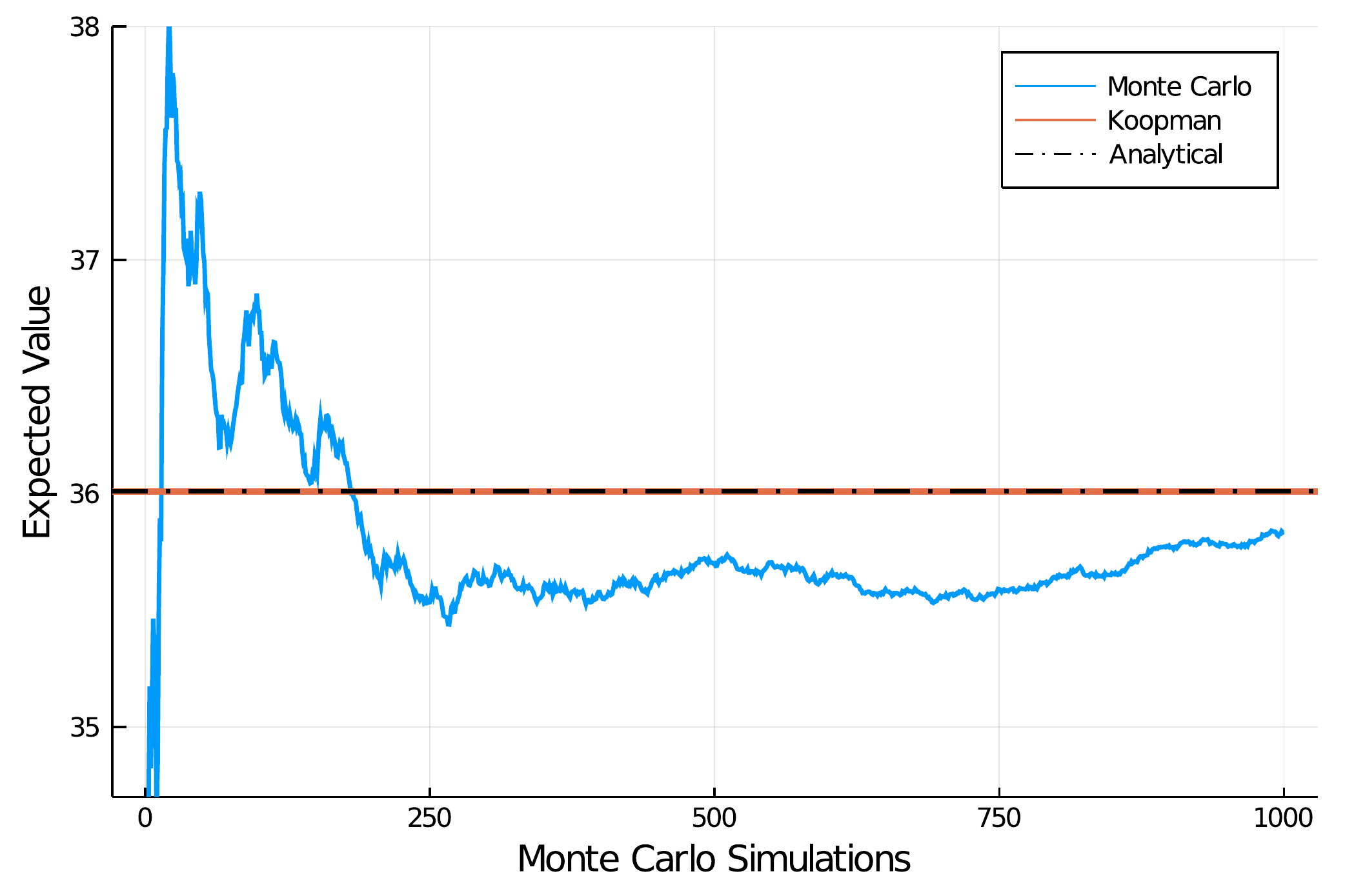}
	\caption{Monte Carlo Expectation Convergence}
	\label{fig:MC_v_Koopman}
\end{figure}

Figure \ref{fig:MC_v_Koopman} shows the convergence of Monte Carlo with respect to the Koopman Expectation solution. After \SI{100000}{} simulations the Monte Carlo solution is still asymptoting to the analytical solution. The absolute error after \SI{100000}{} Monte Carlo simulations is \SI{0.226}{\meter\squared}. Whereas the absolute error resulting from the Koopman Expectation is \SI{2.68e-11}{\meter\squared}.

Next, we leverage this speed-up to optimize $x_0\in\sbr{\SI{-100}{\meter},\SI{0}{\meter}}, \dot{x}_0\in\sbr{\SI{1}{\meter\per\second},\SI{3}{\meter\per\second}}, z_0\in\sbr{\SI{10}{\meter},\SI{50}{\meter}}$ such that $\E{g}$ is minimized. This is achieved using the Method of Moving Asymptotes (MMA) \cite{svanberg_class_2002} gradient-based local optimization algorithm with a relative tolerance stopping criteria on the optimization parameters of \SI{1e-3}{}. Gradients of the expected value are computed using forward mode automatic differentiation. 

\begin{figure}[h]
 	\centering
 	\includegraphics[width=0.65\linewidth]{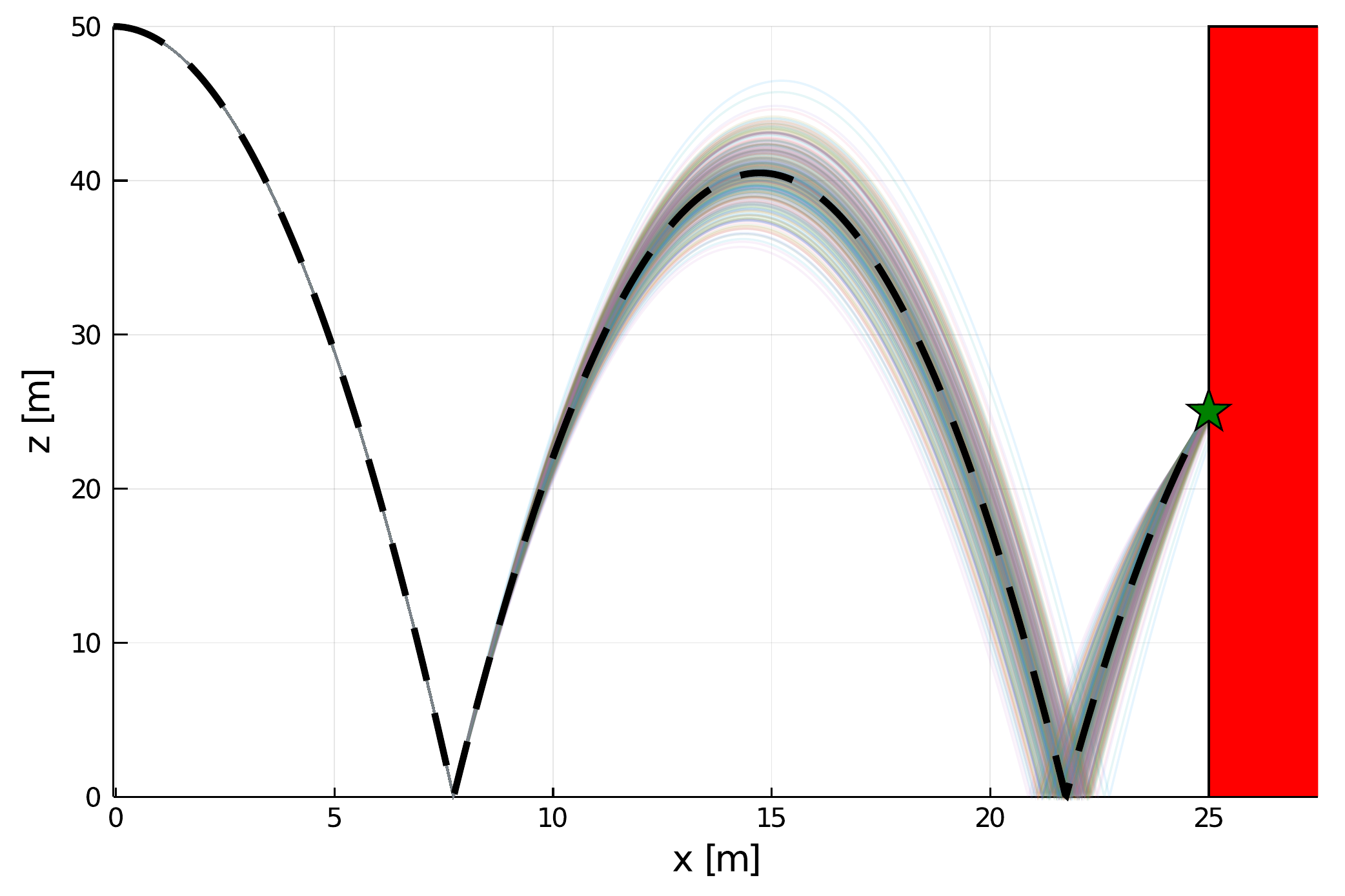}
 	\caption{Bouncing Ball, Minimized Expected Squared Miss Distance}
 	\label{fig:ball_opt}
 \end{figure}

The resulting solution produces an expectation of \SI{8.38e-2}{\meter\squared} in \SI{0.117}{\second} with 26 evaluations of the objective function and gradient in total. By leveraging the Koopman Expectation within the optimization loop, we are able to complete this optimization 17.5x faster than the time required to compute the objective function once using Monte Carlo. Figure \ref{fig:ball_opt} shows the nominal trajectory (dashed line) and 350 random trajectories resulting from the optimized initial conditions.

Lastly, we are interested in computing some higher-order statistics of the observable $g$ resulting from this optimization. In particular, we compute central moments 2-5 using Eq.~\ref{eq:central_moment_expanded} by first defining a vector-valued observable as
\begin{equation}
	\bar{\bs g}\del{\bs x} =\sbr{g\del{\bs x}, g\del{\bs x}^2, ...,g\del{\bs x}^5}^\top \label{eq:obs_vector}
\end{equation}
and then formulate the Koopman Expectation by 
\begin{align}
\Egiven{\bar{\bs g}\del{S\del{\bs x,\alpha}}}{f} = \int_{0.84}^{1}\bar{\bs g}\del{S\del{\bs x, \alpha}}f\del{\alpha}d\alpha
\end{align}

\begin{figure}[h]
	\centering
	\includegraphics[width=0.85\linewidth]{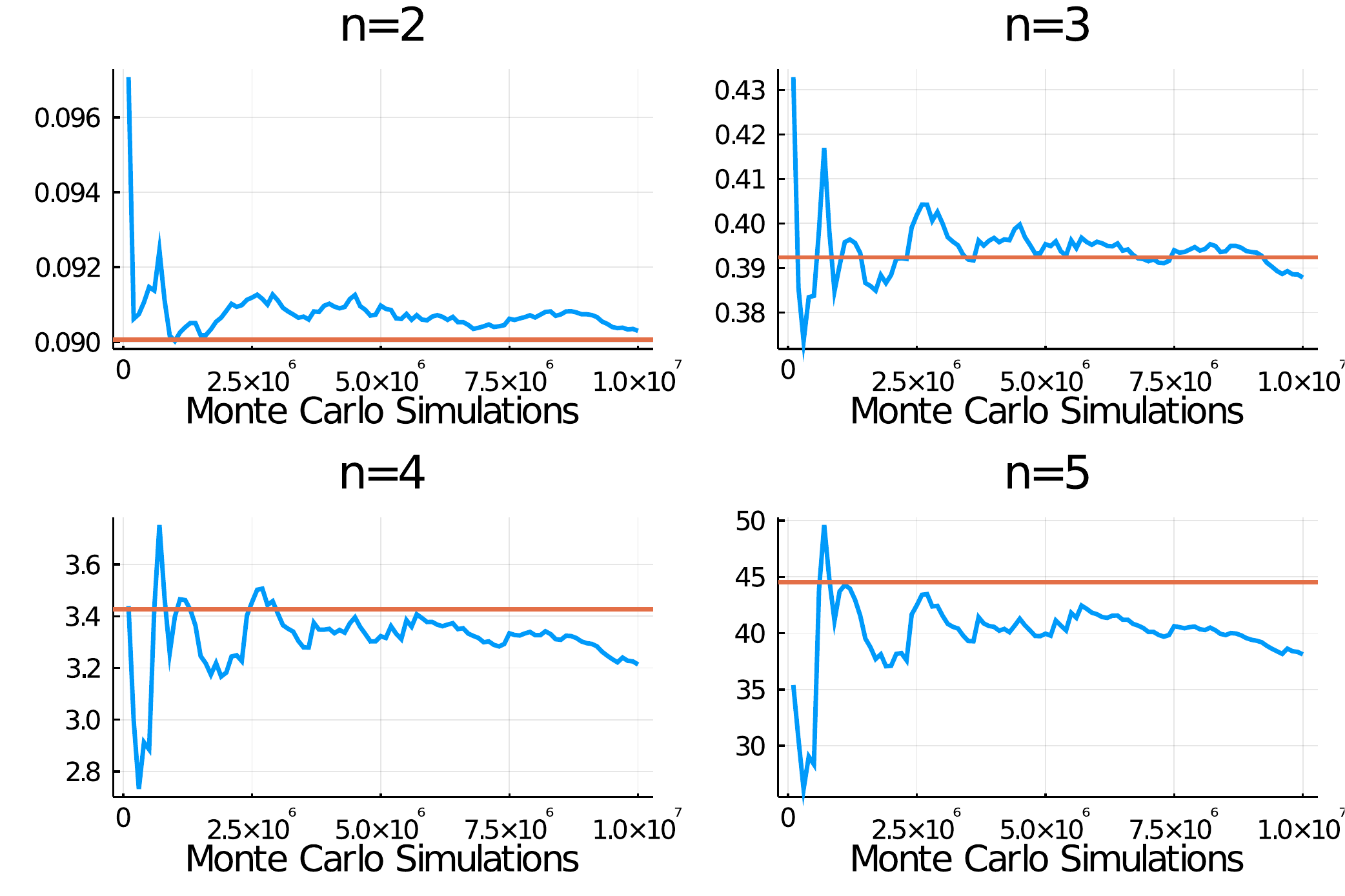}
	\caption{Monte Carlo $n^\mathrm{th}$ Central Moment Convergence}
	\label{fig:ball_moments}
\end{figure}

By leveraging the Koopman Expectation, we are able to compute this vector-valued expectation and then solve Eq.~\ref{eq:central_moment_expanded} for the 4 central moments in $\SI{3.396}{\milli\second}$. This was achieved with only 225 simulations in total. Figure \ref{fig:ball_moments} shows the convergence of Monte Carlo for these moments up to 10M simulations (blue) along with the values computed via the Koopman Expectation (red). Here, 10M Monte Carlo simulations takes \SI{264.5}{\second}. Table \ref{tab:mc_v_koopman_moments} shows the moments computed via the 10M Monte Carlo simulations and the Koopman Expectation. Again, the residual errors computed by the integration method are included for the Koopman Expectation. Compared to 10M Monte Carlo simulations, the Koopman Expectation provides a  77,000x speed-up while realizing tight error bounds on the solution. 


\begin{table}[ht]
	\centering
	\caption{Bouncing Ball Central Moments}
	\label{tab:mc_v_koopman_moments}
	\begin{tabular}[t]{ccc}
		\toprule
		Central Moment & Monte Carlo  & Koopman \\ 
		\midrule
		2 & \SI{9.030e-2}{}& $\SI{9.007e-2}{} \pm \SI{3.878e-5}{}$\\ 
		3 & \SI{3.878e-1}{} & $\SI{3.924e-1}{} \pm \SI{1.776e-3}{}$ \\ 
		4 & 3.214 &  $\SI{3.428}{} \pm \SI{1.536e-3}{}$\\ 
		5 &  38.116&  $\SI{44.536}{} \pm \SI{3.733e-3}{}$\\
		\bottomrule
	\end{tabular}
\end{table}

%
%
%
%
%
%
%



\section{Conclusion}
The need to propagate uncertainty through a dynamical system is prevalent across the scientific and engineering disciplines. Although many solution approaches exist in the literature that can exploit the structure of the dynamics and/or uncertainties to solve such problems efficiently, few naturally generalize to the broadest class of problems involving non-linear hybrid systems with non-Gaussian uncertainties driven by process noise. The de facto standard for such problems is Monte Carlo methods. 

In this work we developed an efficient method for computing the expectation of random variables as propagated through a dynamic system that generalizes to this broadest class of problems. This is achieved by exploiting the action of the Koopman operator without needing an explicit representation of the operator itself. Although this approach only directly applies to analyzing uncertainty through expectations, the careful selection of observables enables the calculation of higher-order statistics. 

When the Koopman Expectation is compared to naive Monte Carlo simulation for a relatively simple hybrid system, we see that the Koopman Expectation provides a significantly more accurate solution while demanding orders of magnitude less computation time. Furthermore, error estimates of the solution are available as a byproduct of the spatial integration method. In many cases, deriving similar error estimates from Monte Carlo simulations is computationally burdensome or even intractable. For applications requiring expectations in some higher-level processing loop, such as in the cases of optimization and control, the impact of these efficiency gains is amplified, enabling the solution of time-sensitive applications. 

Additional research is required to determine how well the Koopman Expectation extends to problems with a high number of uncertain parameters. This is especially true for problems involving process noise, where many uncertain parameters may be required to appropriately approximate the desired noise process.

\bibliography{AFOSR.bib}

\begin{thebibliography}{10}

\bibitem{andronov_theory_2011}
A.~A. Andronov, A.~A. Vitt, and S.~E. Khaikin.
\newblock {\em Theory of {{Oscillators}}}.
\newblock {Dover Publications}, abridged edition edition edition.

\bibitem{arnold1995random}
Ludwig Arnold.
\newblock Random dynamical systems.
\newblock In {\em Dynamical systems}, pages 1--43. Springer, 1995.

\bibitem{brunton_koopman_2016}
Steven~L. Brunton, Bingni~W. Brunton, Joshua~L. Proctor, and J.~Nathan Kutz.
\newblock Koopman {{Invariant Subspaces}} and {{Finite Linear Representations}}
  of {{Nonlinear Dynamical Systems}} for {{Control}}.
\newblock 11(2):e0150171.

\bibitem{budisic_applied_2012}
Marko Budi\v{s}i\'c, Ryan Mohr, and Igor Mezi\'c.
\newblock Applied {{Koopmanism}}.
\newblock 22(4):047510.

\bibitem{burkardt_truncated_2014}
John Burkardt.
\newblock The {{Truncated Normal Distribution}}.

\bibitem{chen_variants_2012}
Kevin~K. Chen, Jonathan~H. Tu, and Clarence~W. Rowley.
\newblock Variants of {{Dynamic Mode Decomposition}}: {{Boundary Condition}},
  {{Koopman}}, and {{Fourier Analyses}}.
\newblock 22(6):887--915.

\bibitem{vcrnjaric2019koopman}
Nelida {\v{C}}rnjari{\'c}-{\v{Z}}ic, Senka Ma{\'c}e{\v{s}}i{\'c}, and Igor
  Mezi{\'c}.
\newblock Koopman operator spectrum for random dynamical systems.
\newblock {\em Journal of Nonlinear Science}, pages 1--50, 2019.

\bibitem{dutta_nonlinear_2015}
Parikshit Dutta, Abhishek Halder, and Raktim Bhattacharya.
\newblock Nonlinear estimation with {{Perron}}-{{Frobenius}} operator and
  {{Karhunen}}-{{Lo\`eve}} expansion.
\newblock 51(4):3210--3225.

\bibitem{genz_remarks_1980}
A.~C. Genz and A.~A. Malik.
\newblock Remarks on algorithm 006: {An} adaptive algorithm for numerical
  integration over an {N}-dimensional rectangular region.
\newblock {\em Journal of Computational and Applied Mathematics},
  6(4):295--302, December 1980.

\bibitem{gleick_chaos_2008}
James Gleick.
\newblock {\em Chaos: {{Making}} a {{New Science}}}.
\newblock {Penguin Books}, anniversary, reprint edition edition.

\bibitem{halmos_legend_1973}
P.~R. Halmos.
\newblock The {{Legend}} of {{John Von Neumann}}.
\newblock 80(4):382--394.

\bibitem{han2017random}
Xiaoying Han and Peter~E Kloeden.
\newblock {\em Random ordinary differential equations and their numerical
  solution}.
\newblock Springer, 2017.

\bibitem{hoogendoorn_uncertainty_2018}
R.~Hoogendoorn, E.~Mooij, and J.~Geul.
\newblock Uncertainty propagation for statistical impact prediction of space
  debris.
\newblock 61(1):167--181.

\bibitem{jones_whither_2001}
Christopher K. R.~T. Jones.
\newblock Whither {{Applied Nonlinear Dynamics}}?
\newblock In Bj\"orn Engquist and Wilfried Schmid, editors, {\em Mathematics
  {{Unlimited}} \textemdash{} 2001 and {{Beyond}}}, pages 631--645. {Springer}.

\bibitem{koopman_dynamical_1932}
B.~O. Koopman.
\newblock Dynamical {{Systems}} of {{Continuous Spectra}}.
\newblock 18(3):255--263.

\bibitem{koopman_hamiltonian_1931}
B.~O. Koopman.
\newblock Hamiltonian {{Systems}} and {{Transformation}} in {{Hilbert Space}}.
\newblock 17(5):315--318.

\bibitem{korda_data-driven_2020}
Milan Korda, Mihai Putinar, and Igor Mezi\'c.
\newblock Data-driven spectral analysis of the {{Koopman}} operator.
\newblock 48(2):599--629.

\bibitem{kutz_dynamic_2016}
J.~Nathan Kutz, Steven~L. Brunton, Bingni~W. Brunton, and Joshua~L. Proctor.
\newblock {\em Dynamic {{Mode Decomposition}}: {{Data}}-{{Driven Modeling}} of
  {{Complex Systems}}}.
\newblock {SIAM-Society for Industrial and Applied Mathematics}.

\bibitem{lasota_chaos_2013}
Andrzej Lasota and Michael~C. Mackey.
\newblock {\em Chaos, {{Fractals}}, and {{Noise}}: {{Stochastic Aspects}} of
  {{Dynamics}}}.
\newblock {Springer Science \& Business Media}.

\bibitem{leonard_probabilistic_2019-1}
Andrew Leonard.
\newblock Probabilistic {{Methods}} for {{Decision Making}} in {{Precision
  Airdrop}}.

\bibitem{leonard_koopman_2019-1}
Andrew Leonard, Jonathan Rogers, and Adam Gerlach.
\newblock Koopman {{Operator Approach}} to {{Airdrop Mission Planning Under
  Uncertainty}}.
\newblock 42(11):2382--2398.

\bibitem{leonard_probabilistic_2020-1}
Andrew Leonard, Jonathan Rogers, and Adam Gerlach.
\newblock Probabilistic {{Release Point Optimization}} for {{Airdrop}} with
  {{Variable Transition Altitude}}.
\newblock In Print:1--11.

\bibitem{meyers_koopman_2019-1}
Joseph~J. Meyers, Andrew~M. Leonard, Jonathan~D. Rogers, and Adam~R. Gerlach.
\newblock Koopman {{Operator Approach}} to {{Optimal Control Selection Under
  Uncertainty}}.
\newblock In {\em 2019 {{American Control Conference}} ({{ACC}})}, pages
  2964--2971.

\bibitem{mezic_spectral_2005-2}
Igor Mezi\'c.
\newblock Spectral {{Properties}} of {{Dynamical Systems}}, {{Model Reduction}}
  and {{Decompositions}}.
\newblock 41(1):309--325.

\bibitem{mezic_comparison_2004}
Igor Mezi\'c and Andrzej Banaszuk.
\newblock Comparison of systems with complex behavior.
\newblock 197(1):101--133.

\bibitem{narasingam_koopman_2019}
Abhinav Narasingam and Joseph Sang-Il Kwon.
\newblock Koopman {{Lyapunov}}-based model predictive control of nonlinear
  chemical process systems.
\newblock 65(11):e16743.

\bibitem{nolte_tangled_2010}
David~D. Nolte.
\newblock The tangled tale of phase space.
\newblock 63(4):33.

\bibitem{goroff_new_1992}
Henri Poincare.
\newblock {\em New {{Methods}} of {{Celestial Mechanics}}}.
\newblock History of {{Modern Physics}} and {{Astronomy}}. {AIP-Press}.

\bibitem{probe_new_2016}
Austin Probe, Tarek~A. Elgohary, and John~L. Junkins.
\newblock A {{New Method}} for {{Space Objects Probability}} of {{Collision}}.
\newblock In {\em {{AIAA}}/{{AAS Astrodynamics Specialist Conference}}}.
  {American Institute of Aeronautics and Astronautics}.

\bibitem{rackauckas_differentialequations.jl_2017}
Christopher Rackauckas and Qing Nie.
\newblock {{DifferentialEquations}}.jl \textendash{} {{A Performant}} and
  {{Feature}}-{{Rich Ecosystem}} for {{Solving Differential Equations}} in
  {{Julia}}.
\newblock 5(1):15.

\bibitem{rackauckas2017adaptive}
Christopher Rackauckas and Qing Nie.
\newblock Adaptive methods for stochastic differential equations via natural
  embeddings and rejection sampling with memory.
\newblock {\em Discrete and continuous dynamical systems. Series B},
  22(7):2731, 2017.

\bibitem{schmid_dynamic_2010}
Peter~J. Schmid.
\newblock Dynamic mode decomposition of numerical and experimental data.
\newblock 656:5--28.

\bibitem{stark1986probability}
Henry Stark and John~W Woods.
\newblock Probability, random processes, and estimation theory for engineers.
\newblock {\em prpe}, 1986.

\bibitem{svanberg_class_2002}
Krister Svanberg.
\newblock A {{Class}} of {{Globally Convergent Optimization Methods Based}} on
  {{Conservative Convex Separable Approximations}}.
\newblock 12(2):555--573.

\bibitem{tsitouras2011runge}
Ch~Tsitouras.
\newblock Runge–kutta pairs of order 5 (4) satisfying only the first column
  simplifying assumption.
\newblock {\em Computers \& Mathematics with Applications}, 62(2):770–775,
  2011.

\bibitem{tu_dynamic_2014}
Jonathan~H. Tu, Clarence~W. Rowley, Dirk~M. Luchtenburg, Steven~L. Brunton, and
  J.~Nathan Kutz.
\newblock On dynamic mode decomposition: {{Theory}} and applications.
\newblock 1(2):391.

\bibitem{weise_global_2009}
Andrea~Yeong Wei\ss{}e.
\newblock Global sensitivity analysis of ordinary differential equations:
  Adaptive density propagation using approximate approximations.
\newblock Accepted: 2018-06-07T22:43:56Z.

\bibitem{williams_datadriven_2015}
Matthew~O. Williams, Ioannis~G. Kevrekidis, and Clarence~W. Rowley.
\newblock A {{Data}}\textendash{{Driven Approximation}} of the {{Koopman
  Operator}}: {{Extending Dynamic Mode Decomposition}}.
\newblock 25(6):1307--1346.

\end{thebibliography}

\appendix
\appendixpage


\section{Bouncing Ball Analytical Solution } \label{ap:ball_analytical}
First, the vertical velocity of the ball after a bounce can be expressed as $\dot{z}^+ = -\alpha \dot{z}^-$. Considering a constant acceleration model and that the initial impact velocity of the ball with the ground after falling $z_0$ from rest is $\dot{z}_1^- = -\sqrt{2gz_0}$, the vertical velocity after the $i^\mathrm{th}$ bounce is
\begin{equation} \label{eq:viplus_vel}
	\dot{z}_i^+ = -\alpha^i \dot{z}_i^- = \alpha^i \sqrt{2gz_0}
\end{equation}
From Eq.~\ref{eq:viplus_vel}, the length of time the ball is airborne after bounce $i$ can be written as 
\begin{equation} \label{eq:ti_bounce}
	t_i = \alpha^i \sqrt{\frac{8 z_0}{g}}
\end{equation}
with the special case of the initial drop from $z_0$, $t_0 = \frac{1}{2} \sqrt{\frac{8 z_0}{g}}$. Therefore, the total at the $n^{\mathrm{th}}$ bounce can be written as
\begin{align}
	t_n = \sum_{i=1}^n t_i + t_0 &= \sum_{i=1}^n \alpha^i \sqrt{\frac{8 z_0}{g}} + \frac{1}{2} \sqrt{\frac{8 z_0}{g}} \nonumber \\
	&= \sum_{i=1}^n \alpha^{(i-1)} \sqrt{\frac{8 z_0}{g}} - \frac{1}{2} \sqrt{\frac{8 z_0}{g}} \nonumber \\
	&= \sqrt{\frac{8 z_0}{g}}  \left( \frac{1-\alpha^n}{1-\alpha} - \frac{1}{2} \right)
\end{align}

By considering the total time for the ball to reach the right hand side wall at $x = x^*$, $T = x^* /\dot{x}_0$, the number of bounces the ball will complete can be computed by solving $\lfloor T - t_n \rfloor = 0$ for $n$,
\begin{align} 
	b(\alpha) &= \frac{\log\left( 1-(1-\alpha)\left(\sqrt{\frac{g}{8z_0}}T - \frac{1}{2} \right) \right)}{\log(\alpha)} \label{eq:nbouncesf}\\
	n &= \lfloor b(\alpha) \rfloor
\end{align}
which is a function of $\alpha$.

Using the number of full bounces, the remainder time and velocity, $t_r = T - t_n$ and $\dot{z}_r =\alpha^n \sqrt{2gz_0}$, respectively, of the ball may be used to compute the impact height, $H$
\begin{equation} \label{eq:yimpact}
H(\alpha) = \dot{z}_r t_r - \frac{1}{2}g t_r^2
\end{equation}
where $n$, $t_r$, and $\dot{z}_r$ are all functions of $\alpha$.

For this work, Eq.~\ref{eq:nbouncesf} is solved for $n = 2$, resulting in the lower-bound requirement $\alpha > 0.8066$. This is satisfied by the truncated Gaussian distribution over $\alpha$. However, any distribution with compact support within $\left(0.8066, 1 \right]$ may be used. With $n = 2$, Eq.~\ref{eq:yimpact} can be expanded to
\begin{equation} \label{eq:himpact_exp}
	H(\alpha) = \alpha^2 \sqrt{2g z_0} \left(T -  \sqrt{\frac{8z_0}{g}}(\alpha + \frac{1}{2}) \right) - \frac{1}{2}g \left(T -  \sqrt{\frac{8z_0}{g}}(\alpha + \frac{1}{2}) \right)^2
\end{equation} 
Further expanding the quadratic and combining like terms, Eq.~\ref{eq:himpact_exp} is cubic with respect to $\alpha$
\begin{align}
	A = \sqrt{2gz_0}, \: B = \sqrt{\frac{8z_0}{g}}, \: C = T - \frac{1}{2}B \nonumber \\
	H(\alpha) = -AB\alpha^3 + (AC - B^2)\alpha^2 + CB\alpha - \frac{1}{2}gC^2
\end{align}
for $\alpha \in \left(0.8066, 1 \right]$.

Considering that the observable $g\del{\alpha}=\del{H\del{\alpha}-z^*}^2$ will generate a $6^\mathrm{th}$ order polynomial, with coefficients $\left\{a_0, \, a_1, \dots, \, a_6\right\}$, the expected value integral of $g$ and $f$ for $\alpha \sim f$, is
\begin{align} 
I &= \int \left( a_0 + a_1\alpha + a_2\alpha^2 + a_3\alpha^3 + a_4\alpha^4 + a_5\alpha^5 + a_6\alpha^6 \right) f(\alpha) \dif \alpha \nonumber \\
&=  a_0 \int f(\alpha) \dif \alpha  + a_1 \int \alpha f(\alpha) \dif \alpha  + a_2 \int \alpha^2 f(\alpha) \dif \alpha + \dots+ a_6 \int \alpha^6 f(\alpha) \dif \alpha  \label{eq:epoly_full}
\end{align}
It can be seen that Eq.~\ref{eq:epoly_full} is a weighted sum of the first 6 raw moments of $f(\alpha)$. These moments are known in closed form, even for the truncated Normal distribution~\cite{burkardt_truncated_2014}. As such, the expected value for the observable $g$ is
\begin{equation}
	\E{g} = \left[a_0, \, a_1, \, a_2, \, a_3, \, a_4, \, a_5, \, a_6 \right] \begin{bmatrix}
	1 \\ m_1 \\ m_2 \\ m_3 \\ m_4 \\ m_5 \\ m_6
	\end{bmatrix}
\end{equation} 
where $m_i$ is the $i^{\mathrm{th}}$ moment of $f$.

The same approach can be extended for higher-order moments of the observable $g$. However, considering that $g^2$ is a 12$^\mathrm{th}$ order polynomial, the coefficients become large creating numerical issues during the multiplication/summation steps.


\end{document}